\documentclass[12pt]{amsart}

\newtheorem{thm}{Theorem}[section]

\newtheorem{lem}[thm]{Lemma}

\newtheorem{cor}[thm]{Corollary}

\newtheorem{prop}[thm]{Proposition}

\theoremstyle{definition}

\newtheorem{note}[thm]{Note}

\theoremstyle{remark}

\newcommand{\R}{\mathbf{R}}

\newcommand{\ol}[1]{{\overline #1}}

\newcommand{\RP}{\mathbf{RP}}

\newcommand{\C}{{\mathbf C}}

\newcommand{\Z}{{\mathbf Z}}

\renewcommand{\S}{\mathbf{S}}

\renewcommand{\(}{\left(}

\renewcommand{\)}{\right)}

\DeclareMathOperator{\conv}{conv}

\DeclareMathOperator{\codim}{codim}

\begin{document}

\title[Totally skew embeddings]{Totally skew embeddings of manifolds}

\author{Mohammad Ghomi}

\address{School of Mathematics, Georgia Institute of Technologies, 
Atlanta, GA 30332}

\email{ghomi@math.gatech.edu}

\urladdr{www.math.gatech.edu/$\sim$ghomi}

\author{Serge Tabachnikov}

\address{Department of Mathematics, Penn State, University Park, PA 16802}

\email{tabachni@math.psu.edu}

\urladdr{www.math.psu.edu/tabachni}



\date{September 2002. Last Typeset \today.}

\thanks{Research of the first author was partially supported by NSF 
Grant DMS-0204190, and CAREER award
DMS-0332333. The second author was  supported in part by NSF 
Grant DMS-0244720 and a BSF grant.}

\begin{abstract}
We obtain bounds on the least dimension of an affine space that can
contain an
$n$-dimensional submanifold without any pairs of parallel or 
intersecting tangent lines at distinct points. This
problem is closely related to the  generalized vector field problem, 
non-singular bilinear maps, and the
immersion problem for  real projective spaces.
\end{abstract}

\maketitle


\section{Introduction}


A pair of lines in   Euclidean space $\R^N$ are said to be 
\emph{skew} if they are not parallel and do not
intersect. Accordingly, we say that a submanifold $M^n$ of $\R^N$ is 
{\it totally skew} if  every pair of
tangent lines of $M$ are skew, unless they are tangent to
$M$ at the same point. The simplest example is  the cubic curve
$$
\R\ni x\longmapsto \(x,x^2,x^3\)\in\R^3.
$$ 
In this paper we construct other examples, while our main focus is 
the following question:
\emph{ Given a manifold $M^n$, what is the smallest dimension 
$N(M^n)$ such that $M^n$ admits a totally
skew embedding in $\R^N$? } For instance, the cubic curve, together 
with the fact that there are no skew
pairs of lines in $\R^2$, easily yields that
$N(\R^1)=3$. In Section \ref{sec:1} of this paper we start our 
investigations by using some classical
convexity arguments, together with Thom's transversality theorem, to 
establish the following basic bounds:

\begin{thm}\label{thm:1}
  For any manifold $M^n$,
$$ 
2n+1\leq N(M^n)\leq 4n+1.
$$
Indeed,  generically any submanifold
$M^n\subset\R^{4n+1}$ is totally skew. Further, if $M^n$ is closed, then
$N(M^n)\geq 2n+2$.
\end{thm}

In particular, there exist no closed totally skew curves in $\R^3$, or, 
in other words, every immersion of the circle
$\S^1$ into the Euclidean 3-space will have a pair of parallel or 
intersecting tangent lines. On the other hand
it is not  difficult to  show by a direct computation that the curve given
by
$$
\C\supset\S^1\ni z\longmapsto \(z,z^2\)\in\C^2\simeq\R^4
$$ 
is totally skew \cite{Gh}. Thus $N(\S^1)=4$. More generally, as we prove
in 
Section
\ref{sec:bilinear}, one may explicitly construct totally skew spheres 
in every dimension:

\begin{thm} \label{thm:2}    
Let
$B: \R^{n+1} \times
\R^{n+1} \to \R^{m}$ be a symmetric nonsingular bilinear map. Then
$$
\R^{n+1}\supset\S^{n}\ni x\longmapsto \big(x,
B(x,x)\big)\in\R^{n+1}\times\R^m
$$ 
is a totally skew embedding. 	In particular,
$$ 
N\(\S^{n}\)\leq n+m(n)+1,
$$ 
where $m(n)$ is the smallest dimension where a symmetric 
nonsingular bilinear map
$\R^{n+1}\times\R^{n+1}\mapsto\R^m$ exists.
\end{thm}

By ``symmetric"  we mean that $B(x,y)=B(y,x)$, and ``nonsingular" 
means that $B(x,y)=0$ only if
$x=0$ or $y=0$. Constructing explicit examples of such bilinear maps in
Section
\ref{subsec:sphere},
  we obtain

\begin{cor}\label{cor:sphere}
$ 
N(\S^n)\leq 3n+2, \quad\text{and}\quad N\(\S^{2\ell-1}\)\leq 3(2\ell-1)+1.
$
\end{cor}

To obtain more  bounds for $N(M^n)$, note that, since $M^n$ is 
locally homeomorphic to
$\R^n$, any totally skew embedding of $M^n$ induces one also of 
$\R^n$. In other words,
$$
  N(M^n) \geq N(n):=N(\R^n).
$$ 
A number of lower bounds for $N(n)$
  are provided by the following result and its corollaries, which are
proved in Sections
\ref{sec:vector} and \ref{sec:cor} respectively.
  Here $\xi_p$ denotes the
\emph{canonical line  bundle} over the real projective space $\RP^N$, and
$r\xi_p$ is the Whitney  sum of $r$ copies of $\xi_p$.

\begin{thm} \label{thm:3}   
If there exists a totally skew disk 
$D^n\subset\R^N$, then
$(N-n)\xi_{n-1}$ admits $n+1$ linearly independent sections. Thus
$$ 
N(n)\geq r(n)+n,
$$ 
where $r(n)$ is the smallest integer such that $r\xi_{n-1}$ admits 
$n+1$ linearly independent sections.
\end{thm}

Finding the number of linearly independent sections of
  $r\xi_p$ is known as the \emph{generalized vector field problem}. It 
has been thoroughly  studied using
various topological methods  \cite{As, Da1, Da2, Da3, D-G-M,  La, 
L-R}. In particular, \cite{L-R} provides a
table of solutions for $1 \leq p < r
\leq 32$. Using these values we immediately obtain:

\begin{cor}\label{cor:table}
  {\footnotesize
$$
\begin{array}{|c|c|c|c|c|c|c|c|c|c|c|c|c|c|c|c|c|c|}\hline n=
&1&2&3&4&5&6&7&8&9&10&11&12&13&14&15&16&17\\ \hline N(n)  \geq
&3&6&7&12&13&14&15&24&25&27&28&31&36&37&38&48&49\\ \hline
\end{array}.
$$ }\qed
\end{cor}

Note that there are significant jumps in the lower bounds in the 
above table whenever $n$ is a power of
$2$. Further this  table suggests that $N(2^p)\geq 3(2^p)$ for all 
$p$. The next Corollary of Theorem
\ref{thm:3}, which is based on the vanishing of Stiefel-Whitney classes of
$(n+q)\xi_{n-1}$ in dimensions $q$ and higher, shows that this is  the case:

\begin{cor} \label{cor:power}  
Suppose that $N(n)=2n+q$, with $q \leq 
n$, and $n+q=2^r+m$ with $0
\leq m < 2^r$.  Then $q >m$, and hence
$2^r-q\leq n \leq 2^r-1.$ In particular,
$$ 
N\(2^\ell\)\geq 3\(2^\ell\),
$$ 
since $q=n$ when $n=2^\ell$.
\end{cor}

The last inequality  in the above corollary together with the bound 
for  spheres in Corollary
\ref{cor:sphere} yield that
$$ 
6\leq N(2)\leq N\(\S^2\)\leq 8.
$$ 
We show in Section \ref{sec:disk} that, 
indeed, $N(2)=6$.   It is reasonable to expect
that, since $\S^2$ is a closed manifold, $N(2)< N(\S^2)$, so one 
might  conjecture that
$N(\S^2)=7$ or
$8$. It would also be interesting to know how often   $N(n)$ is equal to
$2n+1$, its lowest possible value. It follows  from \cite{Da1} that 
this occurs very rarely:

\begin{cor}\label{cor:lowest}
$ 
N(n)\geq 2n+2,  \quad\text{unless}\quad\text{$n=1$, $3$, or $7$}.
$
\end{cor}

Recall that  $N(1)=3$. Thus excluding the case $n=1$ in the 
above corollary is not superfluous;
however, we do not know whether excluding the other two dimensions 
are necessary as well. In other words:
\emph{Does there exist totally skew disks
$D^3\subset\R^7$ or $D^7\subset \R^{15}$}?

The generalized vector field problem is closely related to the 
immersion problem for real projective spaces,
which has been extensively  studied, and a wealth of results is 
known; see \cite{Da4, Ja} for  surveys. 
In particular, one may use the relation between the  immersion 
problem for real projective spaces and
nonsingular maps. Then the strong non-immersion theorem of Davis
\cite{Da3} yields:

\begin{cor}\label{cor:imm}
$\RP^{n-1}$ can be immersed in $\R^{N(n)-n-1}$. In particular,
$$ 
N(2\ell+1)\geq 2\big(3\ell-2d-\alpha(\ell-d)\big)+3,
$$ 
where
$\alpha(m)$ is the number of 1's in the dyadic expansion of
$m$,   and $d$ is the smallest nonnegative integer such that
$\alpha(\ell-d) \leq d+1$.
\end{cor}

In the Appendix we will also present a sharp result for the existence of
\emph{totally skew pairs} of submanifolds in Euclidean space, which 
we will prove using two different
methods. Next we will mention a number of other related results and 
motivations for this paper.

Totally skew submanifolds (TS-embeddings) are generalizations of 
\emph{skew loops}  and
\emph{skew branes}  \cite {G-S,Ta, ghomi:shadow} which are 
submanifolds in  an affine space without parallel
tangent spaces (S-embeddings). Skew loops were first studied by B. 
Segre \cite{segre:global} and have
connections to quadric surfaces. In particular, the first author and 
B. Solomon \cite{G-S} showed that the
absence of skew loops characterizes quadric surfaces of positive 
curvature (including ellipsoids), while the
second author \cite{Ta} ruled out the existence of skew branes on 
quadric hypersurfaces in any dimension.
Note, however, that skew loops are only affinely invariant, whereas 
quadric surfaces are invariant under the
more general class of projective transformations. Thus we were led to 
generalize the concept of skew loops
accordingly. This was the original motivation for our paper.

There is also another type of submanifolds, weaker than totally 
skew,  which are defined by requiring that
there exist no pairs of intersecting tangent lines at distinct points
\cite{Gh}. The main question of this paper may also be asked for such 
submanifolds, which we call
\emph{tangent bundle embeddings} (T-embeddings), as well as for skew 
submanifolds. These questions may
be regarded as generalizations of Whitney's embedding problem
\cite{whitney:2n} who proved that every manifold $M^n$ may be embedded in
$\R^{2n}$.  Of course, TS-embeddings can be viewed as tangent 
bundle embeddings
  into real projective spaces.

  Finally we mention a  somewhat related  notion of 
\emph{$k$-regular} embeddings studied in
\cite {B-R-S, Ch, Co-Ha, Ha1, Ha2, Ha3, Ha-S}.   A  submanifold $M^n
\subset \R^N$ is called
$k$-regular if every $k$-tuple of its distinct points spans  a 
$k$-dimensional vector subspace (there is also an
affine version of this notion). For example, it is proved in
\cite{Co-Ha} that there are no
$k$-regular embeddings of $\R^2$ into $\R^{2k-\alpha(k)-1}$ where, as
before,
$\alpha(k)$ is the number of 1's in the dyadic expansion of
$k$; when $k$ is a power of 2, this result is best possible.


\section{Basic Bounds: Proof of Theorem \ref{thm:1}}\label{sec:1}


\subsection{} First we obtain the lower bounds. 
A pair of affine subspaces $V$,
$W\subset\R^N$ are  \emph{skew} provided that  every pairs of  lines 
$\ell_1\in V$ and
$\ell_2\in W$ are skew.
    Let $AG_n(N)$ denote the (affine Grassmanian) manifold of the
$n$-dimensional affine subspaces of $\R^N$, and $G_{n+1}(N+1)$ be the 
(Grassmanian) manifold of the
$(n+1)$-dimensional (vector) subspaces of
$\R^{N+1}$. There exists a canonical embedding
\begin{equation}\label{eq:i} 
i\colon AG_n(N)\to G_{n+1}(N+1),
\end{equation}  
given by assigning to each point $p\in\R^N$  the  line
$\ell(p)\subset\R^{n+1}$ which passes through the origin and 
$(p,1)$. 
The following observation is
immediate:

\begin{lem}\label{lem:A&G} 
Two affine subspaces $V$, $W\subset\R^N$ 
are  skew, if, and only if, the
corresponding subspaces $i(V)$, $i(W)\subset\R^{N+1}$ have no 
nontrivial intersection, that is,
$i(V)\cap i(W)=0$. In particular, if $V^n$, $W^m\subset\R^N$ are  skew
affine  subspaces, then
$N\geq m+n+1$.\qed
\end{lem}

In particular, no pairs of distinct tangent spaces to a submanifold
$M^n\subset\R^{2n}$ are skew. So we conclude that  $N(M^n)\geq 2n+1$. 
The next observation improves
this bound in the case where $M$ is \emph{closed}, i.e., compact, 
connected, and without boundary.

\begin{prop}\label{prop:support} 
Let $M^n$ be a compact manifold, and 
$f\colon M\to\R^{n+k}$ be a
continuous map, where $k\geq 2$. Then there exists a pair of distinct
points
$p$,
$q\in M$ and a hyperplane $H\subset\R^{n+k}$ such that $f(p)$, 
$f(q)\in H$, while $f(M)$ lies entirely on
one side of $H$.
\end{prop}

\begin{proof} 
Let $C$ be the convex hull of  $f(M)$. If $C$ 
has no interior points, i.e., it lies in a
hyperplane, then we are done. Otherwise $\partial C$ is homeomorphic 
to $\S^m$, where 
$m=n+k-1> n$. Let
$X:= f^{-1}(\partial C)$. If $f|X$ is not one-to-one, then we are 
done, because through each point of
$\partial C$ there passes a support plane. Suppose then that $f$ is 
one-to-one on
$X$. Then since $X$ is compact, and $\partial C$ is  
Hausdorff, it follows that $f(X)$ is homeomorphic
to
$X$. In particular, $f(X)\neq \partial C$, because $\partial 
C\simeq\S^m$, and  $\dim(M) < m$. So
there exists a point $r\in\partial C$ such that $r\not\in f(M)$. It 
follows then from Caratheodory's
  theorem \cite{Sch} that $r$ must lie in the interior of a simplex
$S$ whose vertices are points of  $f(M)$. Let $H$ be a support 
plane of $C$ at
$r$. Then $H$ also supports $S$, and therefore $S$ must lie in $H$ 
because $H$ contains an interior point of
$S$. But if $S$ lies in $H$ so do its vertices. So we obtain at least 
two points of $M$ which lie in
$H$.
\end{proof}

Suppose that $M^n\subset\R^{2n+1}$ is a closed submanifold. Then, letting
$f$ be the inclusion map, Proposition \ref{prop:support} implies that 
there exists a support hyperplane
$H$ of
$M$ which intersects $M$ at distinct points $p$ and $q$. Since $M$ 
has no boundary,  $H$ must be tangent
to $M$ at these points. Thus
$M$ has a pair of distinct tangent spaces $T_p M$ and $T_q N$ which both lie
in
$H\simeq\R^{2n}$. So it follows  from Lemma
\ref{lem:A&G} that $T_p M$ and $T_q N$ must contain a pair of parallel 
or intersecting tangent lines. We
conclude, therefore, that $N(M^n)\geq 2n+2$, when $M$ is closed.

\subsection{} 
Now we derive the generic upper bound for $N(M^n)$.
Since every manifold $M^n$ may be easily embedded in $\R^{2n+1}$ 
\cite{H}, the following result yields that
$N(M^n)\leq 4n+1$.

\begin{prop} 
Any $C^1$-immersed submanifold $M^n\subset\R^{4n+1}$ 
becomes totally skew after an
arbitrary small perturbation.
\end{prop}

\begin{proof} 
For any $C^1$ map $f\colon M\to\R^N$, and pairs of 
points $p$, $q\in M$, the corresponding
$1$-multijet of $f$ is given by
$$ 
j^1_{[2]}f(p,q):=\big(\,p, q, f(p), f(q), df_p, df_q\,\big),
$$ 
where $d$ is the differential map. The set of all such jets, ranging over
  mappings $f$ and pairs of points of $M$, is denoted by
$J^1_{[2]}(M,\R^N)$. Thus we may write
$$ 
M\times M\ni(p,q)\overset{j^1_{[2]}f}{\longmapsto} 
j^1_{[2]}f(p,q)\in J^1_{[2]}
\left(M,
\R^N\right).
$$

  Let $A\subset J^1_{[2]} (M, \R^N)$ be  the set of those jets  of
immersions for which
$df_p(T_pM)$  and
$df_q(T_qM)$  are not  skew. In other words, by Lemma \ref{lem:A&G},
$$ 
A:=\Big\{\,j^1_{[2]}f(p,q) \;\big|\; i\big(df_p(T_pM)\big)
\cap i\big(df_q(T_qM)\big)\neq 0\,\Big\}.
$$
  Then $f(M)$ is totally skew if and only if the jet extension of $f$ 
is completely  disjoint from $A$:
$$
 j^1_{[2]} f(M\times M)\cap A=\emptyset.
$$

By Thom's  transversality theorem \cite{G-G,eliashberg}, if
$A$ is a stratified subset, then, after a perturbation of $f$, we may 
assume that
$j^1_{[2]} f(M\times M)$ is transversal to $A$. Thus once we check 
that $A$ is stratified,  then to establish
the above equality, for generic
$f$, it would suffice to show that
$$
\dim \Big(j^1_{[2]} f(M\times M)\Big)+\dim(A)<\dim\Big(J^1_{[2]} \left(M,
\R^N\right)\Big).
$$ 
Since $\dim(M)=n$,   this is equivalent to
$$
\codim(A):=\dim\Big(J^1_{[2]} \left(M,\R^N\right)\Big)-\dim(A)>2n.
$$

To establish the last inequality,  note first that $J^1 (M,\R^N)$ has 
the  local trivialization
$\R^n \times L(\R^n, \R^N)$ where $L(\R^n, \R^N)$ is the space of 
affine injections $\R^n \to \R^N$. One
also has a fibering $\pi: L(\R^n, \R^N) \to AG_n(N)$,  assigning  to an 
affine map its image. Identifying
$df_p$ with $i(df_p(T_pM))$, we obtain a local projection
$$ 
\pi \times \pi: J^1_{[2]} (M,\R^N) \to G_{n+1}(N+1)\times G_{n+1}(N+1).
$$ 
Let $B \subset G_{n+1}(N+1)\times G_{n+1}(N+1)$ consist of pairs 
of  subspaces in $\R^{N+1}$ with a
nontrivial intersection. Then $A= (\pi \times \pi)^{-1} (B)$, and it follows 
that $\codim(A)$ in
$J^1_{[2]} (M,\R^N)$ equals $\codim(B)$ in
$G_{n+1}(N+1)\times G_{n+1}(N+1)$.

Fix an element of $G_{n+1}(N+1)$, say
  $\R^{n+1}$. Then almost any   $V\in G_{n+1}(N+1)$ is the graph of a 
linear transformation
$\phi\colon\R^{n+1}\to \R^{(N+1)-(n+1)}$. Note that $V_1\cap V_2\neq 
0$, if and only if, ${\rm
rank}(\phi_1 - \phi_2)<n+1$. Thus $B$ is locally identified with the 
set of pairs $(\phi_1,
\phi_2)$ of operators $\R^{n+1}\to \R^{N-n}$ with ${\rm rank}(\phi_1 
- \phi_2)<n+1$. Let
$C$ be the  set of linear maps $\R^{n+1}\to
\R^{N-n}$ whose rank is less than $n+1$. This is  an algebraic 
variety and is therefore stratified. Further,  by
the ``corank formula" (see, e.g, \cite{G-G}) , the  codimension of $C$ is
$((N-n)-n)((n+1)-n)= N-2n$. One has an obvious linear projection $B 
\to C$ sending $(\phi_1,
\phi_2)$ to $\phi_1 - \phi_2$. This projection induces a 
stratification of $B$, and one concludes  that
$\codim(B)$ in $G_{n+1}(N+1)\times G_{n+1}(N+1)$ is also $N-2n$. Thus 
it follows that
$\codim(A)=N-2n$. In particular, if $N\geq 4n+1$, then
$\codim(A)\geq2n+1>2n$, as desired.
\end{proof}


\section{Bilinear Maps: \\Proof of Theorem \ref{thm:2} and Its
Corollary}\label{sec:bilinear}


\subsection{}  
Let $M\subset\R^N$ be a submanifold and
$p\in M$. Then by
$T_p M$ we denote the tangent space of $M$ at $p$, while $\ol T_p M$ 
denotes the translation of
$T_p M$ to the origin, that is
$$ 
T_pM=p+\ol T_p M.
$$

\begin{lem}\label{lem:B'} Let $B\colon\R^{n+1}\times\R^{n+1}\to\R^m$ 
be a symmetric  bilinear map,
$X:=\big(x, B(x,x)\big)$, and
$M:=\{X\mid x\in\S^n\}$. Then,
$$
\ol T_X M=
\Big\{\,\big(v, 2B(x,v)\big) \;\big|\; v\in 
\R^{n+1},\;\text{and}\;\langle v, x\rangle=0\,\Big\}
$$ 
for every $X\in M$.
\end{lem}

\begin{proof} 
First note that for any $C^1$ curve $c\colon(-\epsilon,
\epsilon)\to\R^{n+1}$,
\begin{eqnarray*}
\frac{d}{dt} B\big(c(t), c(t)\big)\Big|_{t=0}&=&
\lim_{t\to 0}\frac{B\big(c(t),c(t)\big)-B\big(c(0),c(0)\big)}{t}\\ 
&=& \lim_{t\to
0}B\left(\frac{c(t)-c(0)}{t},c(t)\right)+\lim_{t\to 
0}B\left(c(0),\frac{c(t)-c(0)}{t}\right)\\
&=&2B\big(c'(0),c(0)\big).
\end{eqnarray*}
 Now suppose that 
$\|c(t)\|=1$, and define
$\gamma\colon (-\epsilon,\epsilon)\to M$ by
$$
\gamma(t):=\Big(c(t), B\big(c(t),c(t)\big)\Big).
$$ 
Suppose $c(0)=x$, and set $v:=c'(0)$. Then $\gamma(0)=X$,
$$
 T_X M\ni\gamma'(0)=\Big(c'(0), 2B\big(c'(0),c(0)\big)\Big)=\big(v, 
2B(v,x)\big),
$$ 
and it remains to note that, since $\|c(t)\|=1$, $\langle v, x\rangle=0$.
\end{proof}

\subsection{} 
Suppose (towards a contradiction) that $M$ is not 
totally skew. Then, for some distinct
$X$, $Y\in M$, one  of the following two conditions holds:

\smallskip

$(i)$ $T_X M$ and $T_Y M$ contain parallel lines,

\smallskip

$(ii)$ $T_X M \cap T_Y M\neq\emptyset$.

\smallskip

\noindent  The above conditions in turn respectively imply that

\smallskip

$(i)'$ $\ol T_X M \cap \ol T_Y M \neq 0$,

\smallskip

$(ii)'$ $Y-X \in \ol T_X M + \ol T_Y M$.

\smallskip

\noindent The implication $(i)\implies (i)'$ is clear. To see
$(ii)\implies (ii)'$, recall that $T_X M=X+\ol T_X M$, and
$T_Y M=Y+\ol T_Y M$. Thus $(ii)$ implies that $X+v=Y+w$ for some
$v\in \ol T_X M$ and $w\in \ol T_Y M$. Consequently,
$Y-X=v-w\in\ol T_X M + \ol T_Y M$.

First assume that $(i)'$ holds. Recall that, by definition of $M$,
\begin{equation}\label{eq:XandY} 
X=\big(x,B(x,x)\big), \quad 
Y=\big(y, B(y,y)\big),
\end{equation}  for some distinct $x,y\in \S^n$.
  Thus, by Lemma \ref{lem:B'}, there exist
\emph{non-zero} vectors
$u,v
\in
\R^n$ such that
$$
\big(v, 2B(x,v)\big) = \big(u, 2B(y,u)\big).
$$ 
So
$u=v$, and $B(x,v)=B(y,u)$. These imply that $B(y-x,v)=0$. Since $B$ 
is non-singular and $x
\neq y$, we conclude that $v=0$, which is a  contradiction.

Next assume that $(ii)'$ holds. Then,  by \eqref{eq:XandY} and Lemma
\ref{lem:B'}, there exist $u,v \in \R^n$ such that
$\langle v, x\rangle = 0$, $\langle u,y\rangle = 0$,  and
\begin{eqnarray} 
y-x=u+v,\ \ B(y,y)-B(x,x)= 2B(x,v)+2B(y,u). \label{cond}
\end{eqnarray} 
Substituting $y=x+u+v$ into the second equation in 
(\ref{cond}) and  collecting terms we
obtain:
$B(u,u)=B(v,v)$. It follows that $B(u+v,u-v)=0$. Thus, since $B$ is 
non-singular, either
$u=-v$ or
$u=v$. If $u+v=0$, then, by the first equation in (\ref{cond}), 
$y=x$,  which contradicts the assumption that
$x$ and $y$ are distinct. If $u=v$ then $v$ is  orthogonal to $x$ and 
$y$, and, by the first equation in
(\ref{cond}), $y-x=2v$.  Taking  the inner-product of both sides of 
this equation with $v$ yields  $0=2
|v|^2$. Thus $v=0$, which is again a contradiction. So we conclude 
that $M$ is totally skew.

\subsection{Proof of Corollary \ref{cor:sphere}}\label{subsec:sphere} Note
that
\begin{equation}\label{eq:B} \Big((x_0,\ldots ,x_n),(y_0,\ldots, 
y_n)\Big)\mapsto\(\sum_{i+j=0} x_i\,
y_j,\ldots,\sum_{i+j=2n} x_i\, y_j\)
\end{equation} gives a nonsingular bilinear map
$\R^{n+1}\times\R^{n+1}\mapsto\R^{2n+1}$. Thus, by  Theorem 
\ref{thm:2}, we have
$$ 
N(\S^n)\leq n+(2n+1)+1=3n+2,
$$
 as had been claimed. Further, \eqref{eq:B} may also be viewed as a 
complex mapping
$\C^{n+1}\times \C^{n+1}\mapsto\C^{2n+1}$, which in turn yields a 
real nonsingular bilinear map 
$\R^{2\ell}\times\R^{2\ell}\mapsto\R^{4\ell-2}$.  Consequently, 
$$
 N\(\S^{2\ell-1}\)\leq (2\ell-1)+(4\ell-2)+1=3(2\ell-1)+1.
$$
\qed

\begin{note}
The  assumption in Theorem \ref{thm:2} that $B$ be symmetric is essential. For
instance, the (nonsymmetric) nonsingular bilinear map $\R^4\mapsto\R^4$ given by
 quaternion multiplication fails to produce a totally skew embedding of $\S^3$ in
$\R^8$, in contrast to the totally skew embedding of $\S^1$ in $\R^4$ given by
complex multiplication.
\end{note}

\begin{note} 
We do not know whether there exist symmetric  bilinear 
non-singular maps
$\R^{n+1} \times \R^{n+1} \mapsto \R^{n+q}$ with $q < n$. 
Construction of such maps would improve the
upper bounds for $N(\S^n)$ obtained in this section. The existence of 
a symmetric  bilinear non-singular map
$\R^n
\times \R^n \mapsto \R^{n+q}$ also implies that $\RP^{n-1}$ can be embedded
in
$\R^{n+q-1}$, see
\cite{Ja}. Thus totally skew submanifolds are  related to
\emph{embeddings} of real projective spaces (in addition to the 
\emph{immersion} problem which will be
discussed in Section 5.3).
\end{note}


\section{The  Vector Field Problem:\\Proof of Theorem
\ref{thm:3}}\label{sec:vector}


  Let $N:=N(n)$. Then there exists a totally skew disc $D^n
\subset \R^N\textsf{}$. Define the (Gauss) map
$$
 D\ni x\overset{f}{\longmapsto} i(T_x D)\in G_{n+1}(N+1),
$$
  where
$i$ is as in
\eqref{eq:i}. By Lemma \ref{lem:A&G}, since $D$ is totally skew,
$$
 f(x) \cap f(y) =0,\quad \text{for all $x \neq y$ in $D$.}
$$ 
 Replacing $D$ by a smaller subset, we may assume that $D$ is the 
graph of a  mapping
$\phi: U^n \to \R^{N-n}$, where $U
\subset \R^n$ is an open disc centered at the origin. Then, for every 
$x \in D$,  $f(x)$  is the graph of a
uniquely determined linear map $A(x): \R^{n+1} \to \R^{N-n}$.

Note that $f(x)\cap f(y)\neq 0$,  if, and only if,  $A(x)-A(y)$ has 
non-zero kernel. Thus
$A(x)-A(y): \R^{n+1} \to \R^{N-n}$ is a linear injection for all $x
\neq y$ in $D$. Such mappings may be identified with the (Stiefel)
manifold
$V_{n+1}(N-n)$ of
$(N-n)\times (n+1)$ matrices with  linearly independent 
columns. Hence we obtain a
  map
$$
 (D \times D - \Delta_D)\ni (x,y)\overset{F}{\longmapsto} 
A(x)-A(y)\in V_{n+1}(N-n),
$$
 where $\Delta_D$ denotes the diagonal elements of $D\times D$. 
Note that $F$ is antisymmetric, that is,
\begin{equation}\label{eq:F} 
F(y,x)=-F(x,y).
\end{equation}
 After a dilation, we may assume that $\S^{n-1}\subset 
U$. Thus we may define a mapping
$$
\S^{n-1}\ni u\overset{\Phi}{\longmapsto}
\Big(\big(u,\phi(u)\big),\big(-u,\phi(-u)\big)\Big)\in D
\times D - \Delta_D.
$$ 
Now if we set $G:=F\circ\Phi$, we obtain a mapping
$$
  \S^{n-1} \overset{G}{\longmapsto} V_{n+1}(N-n).
$$ 
Note that $\Phi(-u)$ is the transpose of $\Phi(u)$, that is, changing $u$ to $-u$
switches the two components of $\Phi(u)$. This, together with \eqref{eq:F}, 
yields that
$G$ is odd:
\begin{equation}\label{eq:G} 
G(-u)=-G(u).
\end{equation}

  Next,  let $\{e_1,\dots, e_{n+1}\}$ be the standard basis for
$\R^{n+1}$ and, for
$1\leq i\leq n+1$, define the mappings $\sigma_i$ by
$$
\RP^{n-1}:=\S^{n-1}/\Z_2\ni 
[u]\overset{\sigma_i}\longmapsto\big[\big(u, G(u)\cdot 
e_i\big)\big]\in
\Big(\S^{n-1} \times\R^{N-n}\Big)\Big/\Z_2,
$$ 
where  $[*]:=\{*,-*\}$.  Note that, using \eqref{eq:G}, we have
$$
\sigma_i\big([-u]\big)=
\big[\big(-u, G(-u)\cdot e_i\big)\big]=
\big[-\big(u, G(u)\cdot e_i\big)\big]=
\big[\big(u, G(u)\cdot e_i\big)\big]=
\sigma_i\big([u]\big).
$$ 
Thus $\sigma_i$ is well-defined. Finally, it remains to note that
$$
\Big(\S^{n-1} \times\R^{N-n}\Big)\Big/\Z_2\simeq (N-n)\xi_{n-1}.
$$  
So  $\sigma_i$ is a section of $(N-n)\xi_{n-1}$ over $\RP^{n-1}$. 
Further, since
$\{G(u)\cdot e_i\}$ is the set of columns of $G(u)$, it is linearly 
independent by definition, and it follows that
$\{\sigma_i\}$ is linearly independent. So we conclude that
$(N-n)\xi_{n-1}$ admits $n+1$ linearly independents sections.


\section{A totally skew disk in $\R^6$}\label{sec:disk}


As we mentioned in the introduction, Theorem \ref{thm:3}, which we
just proved, immediatly yields Corollary
\ref{cor:table}, and, in particular, the bound $N(2)\geq 6$.
In other words, the smallest possible dimensions of an affine space which
contains a totally skew disk cannot be less than $6$. On the other hand:

\begin{prop} \label{prop:dim2} The complex cubic curve given by
$$
\R^2\simeq\C\ni z\longmapsto \(z,z^2,z^3\)\in\C^3\simeq\R^6
$$
 is a totally skew embedding. In particular, $N(2)\leq 6$.
\end{prop}

\begin{proof} 
First we show that the complex cubic curve is totally 
skew in the complex sense, that is,  it has
neither intersecting nor parallel complex tangent lines. Indeed, the 
tangent vector at point
$P(z):=(z,z^2,z^3)$ is $(1,2z,3z^2)$. If two such vectors are 
parallel at points $P(z)$ and
$P(z_1)$ then $z=z_1$. If the tangent lines at $P(z)$ and
$P(z_1)$ intersect then
$$
\(z-z_1, z^2-z_1^2, z^3-z_1^3\)=u\(1,2z,3z^2\)+v\(1,2z_1,3z_1^2\)
$$
 for 
some $u,v \in \C$ (cf. the proof
of Theorem \ref{thm:2} above). Equating the first components, we have
$z-z_1=u+v$; equating the second yields $z^2-z_1^2=2uz+2vz_1$. Hence 
$(u+v)(z+z_1)= 2uz+2vz_1$,
and therefore $(u-v)(z-z_1)=0$. It follows that either $z=z_1$ or 
$u=v$. In the later case equate the third
components to obtain $z^3-z_1^3 = 3u(z^2 + z_1^2)$, and since
$z-z_1=2u$, one has $2u(z^2 + z z_1 + z_1^2) = 3u(z^2 + z_1^2)$. It 
follows that
$u(z-z_1)^2=0$. Therefore either $z=z_1$ or $u=0$, which again 
implies that $z=z_1$.

Suppose now that the 2-dimensional disc is not totally skew. Then 
either the real tangent lines at some
distinct  points $P(z)$ and  $P(z_1)$ are parallel or they intersect. 
In the former case, let $\xi$ and $\xi_1$
be  parallel tangent vectors. Then the vectors $J\xi$ and $J\xi_1$ 
are also parallel where $J$ is the complex
structure in $\C^3$, the operator of multiplication by $\sqrt{-1}$. 
It follows that the complex tangent  lines
at $P(z)$ and  $P(z_1)$ are parallel.  On the other hand, if the real 
tangent lines at $P(z)$ and
$P(z_1)$ intersect, then so do the complex tangent lines as well. In 
both cases we obtain a  contradiction.
\end{proof}

So the above proposition together with Corollary \ref{cor:table} yield:

\begin{cor}
$N(2)=6$.\qed
\end{cor}


\section{Proof of the Other Corollaries of Theorem
\ref{thm:3}}\label{sec:cor}


\subsection{Proof of Corollary \ref{cor:power}}
  The total Stiefel-Whitney class of  $(n+q)\xi_{n-1}$ is $(1+x)^{n+q}$,
where
$x$ is the generator of the homology group
$H^1(\RP^{n-1}; \Z_2)$. Since the  Stiefel-Whitney classes vanish in 
dimensions $q$ and higher,
\begin{eqnarray}
   {{n+q} \choose i} =0\ ({\rm mod}\ 2)\quad {\rm for}\quad i = 
q,\dots, n-1. \label{bin}
\end{eqnarray} 
By assumption,  $2^{r+1} > n+q  =2^r + m$ and $n \geq 
q$. These inequalities imply that
$2^{r+1}>2q$, and hence $2^r>q$. It follows that
$n-m=2^r-q>0$. Thus $m \leq n-1$.

Suppose $q \leq m$. Then $m \in \{q,\dots, n-1\}$. It is well  known 
that the maximal power of $2$ dividing
$k!$ equals $k-\alpha(k)$. It follows that
$$
{{2^r+m} \choose m} =1\ ({\rm mod}\ 2) 
$$ 
which contradicts 
(\ref{bin}). So we conclude that that
$q>m$.\qed

\subsection{\bf Proof of Corollary \ref{cor:lowest}}
   In \cite{Da1},  Davis  estimated the maximal number of linearly 
independent sections of the bundles
$k\xi_p$. In  particular, if ${{k-1}
\choose p}$ is odd then $k\xi_p$ has at most
$$
s:=k-p+2\nu(k)+\epsilon\big(\nu(k),p\big)
$$ 
independent sections; here 
$\nu(k)$ is the greatest power of
$2$  dividing $k$, and $\epsilon$ depends on the mod $4$ values of 
its arguments with
$$
\epsilon(0,2)=\epsilon(3,2)=0,\ \epsilon(1,2)=\epsilon(2,2)=-2.
$$
 Assume
that
$N=2n+1$, that is, $q=1$. By Corollary \ref{cor:power},
$n=2^r-1$. So, by Theorem
\ref{thm:3}, $2^r \xi_{2^r-2}$ is a trivial bundle. Apply the theorem of Davis 
 with $k=2^r$ and
$p=2^r-2$. Then $\nu(k)=r, p=2\ {\rm mod}\ 4$ and $s=2+2r$ for
$r=0,3\ {\rm mod}\ 4$ and $s=2r$ for $r=1,2\ {\rm mod}\ 4$. Since 
$2^r \xi_{2^r-2}$  has
$2^r$  sections, one has: $s \geq 2^r$, and therefore $2+2r \geq 
2^r$, which is clearly impossible unless $r
\leq 3$.\qed

\subsection{Proof of Corollary \ref{cor:imm}} By Theorem \ref{thm:3}, 
$(N(n)-n)\xi_{n-1}$ admits
$n+1$ linearly independent sections. As we showed in Section 4.1, 
this implies that there exists a
$\Z_2$-equivariant map
$G: \S^{n-1} \to V_{n+1}(N-n)$. Since $V_{n+1}(N-n)$ may be 
identified with the space of  linear injections
$\R^{n+1} \to
\R^{N-n}$, $G$ induces a mapping $\S^{n-1} \times \R^{n+1} \mapsto
\R^{N-n}$  and, by homogeneous extension, a {\it non-singular map}
$$ 
g: \R^n \times \R^{n+1} \to \R^{N-n}.
$$ 
That is, $g(x,y)=0$ only if $x=0$  or $y=0$. Further,  $g$ is 
homogeneous of degree 1 in each variable:
$g(sx,ty)=st\  g(x,y)$ for all $s,t \in
\R$. Equivalently, one  obtains an {\it axial map}
$$ 
{\bar g}: \RP^{n-1} \times \RP^n \to \RP^{N-n-1},
$$ 
where  ``axial" means that the restriction of ${\bar g}$ on  each 
factor is homotopic to the inclusion. By
restriction, one also has an axial map
$$ 
h: \RP^{n-1} \times \RP^{n-1} \to \RP^{N-n-1}.
$$ 
The existence of such a map is equivalent to the existence of an 
immersion $\RP^{n-1}
\to
\R^{N-n-1}$, see \cite{Da4, Ja}. The stated lower bound for $N(n)$ 
now follows from the immersion theorem
of Davis  \cite{Da3}.\qed

\begin{note} The relation between totally skew  submanifolds and 
immersions of projective spaces is rather
unexpected. Another,  seemingly unrelated problem which turned out to 
be equivalent to the immersion
problem for
$\RP^n$ concerns the topological complexity of the motion planning 
problem on $\RP^n$
\cite{F-T-Y}.
\end{note}

\section*{Appendix: Totally Skew Pairs}
We say that submanifolds $M_1$, $M_2 \subset \R^N$  are {\it  a 
totally skew pair}, if, for every $x \in
M_1,\ y \in M_2$, the tangent spaces $T_{x} M_1$ and
$T_y M_2$  are skew. It follows immediately from Lemma \ref{lem:A&G} that
$N\geq n_1+n_2+1$. The following result gives  an improvement 
of this bound for closed submanifolds.

\begin{thm}\label{thm:pair} 
If $M_1^{n_1}$, $M_2^{n_2}\subset \R^N$ 
are  a totally skew
\emph{pair} of
\emph{closed} submanifolds, then
$N\geq n_1+n_2+2$.
\end{thm}

This lower bound is sharp, because there exists a totally skew pair in
$\R^{n_1+n_2+2}= \R^{n_1+1} \times \R^{n_2+1}$, given by
$$
  M_1 := \S^{n_1}\times\{0\}, \quad \text{and} \quad M_2:=\{0\} \times
\S^{n_2} .
$$  
We offer two proofs for Theorem \ref{thm:pair}.

\begin{proof}[Proof 1] Since $M_1\cap M_2=\emptyset$,  the (Gauss) map
$$  
M_1 \times M_2 \ni (x,y)\overset{f}{\longmapsto}
\frac{(x-y)}{\|x-y\|}
\in \S^{N-1}
$$ 
is well-defined. Recall that the differential of $f$, evaluated at 
a tangent vector $(v,w)$ in $\ol
T_{(x_0,y_0)}(M_1\times M_2)$,   is given by
$$ 
df_{(x_0,y_0)}\big((v,w)\big):=\frac{d}{dt} 
f\big(x(t),y(t)\big)\Big|_{t=0}\,,
$$
 where
$t\mapsto(x(t), y(t))$ is a curve in $M_1\times M_2$  with
$(x(0),y(0))=(x_0,y_0)$ and
$(x'(0),y'(0))=(v,w)$. Further note that the definition of $f$ yields
$$
\frac{d}{dt} f(x,y)=\frac{(x'-y')\|x-y\|-(x-y)\langle 
x'-y',x-y\rangle\big/\|x-y\|}{\|x-y\|^2}.
$$ 
Now suppose that $df_{(x_0,y_0)}((v,w))=0$. Then the numerator of 
the above fraction vanishes at
$t=0$, and we obtain
$$ 
v-w=k(x_0-y_0),
$$ 
for some scalar $k$. If $k\neq 0$, then
$
  x_0- {v}/{k}=y_0-w/k,
$
  and hence
$ T_{x_0}M_1$ intersects $T_{y_0}M$ which is a contradiction. So 
$k=0$, which yields
$v=w$. Since $M_1$ and
$M_2$ are a totally skew pair, this can happen only when 
$(v,w)=(0,0)$. So we conclude that
$df$ is nonsingular everywhere. In particular, $N-1\geq n_1+n_2$. 
Thus to complete the proof it remains to
show that $N-1\neq n_1+n_2$.

Suppose that $N-1=n_1+n_2$. Then, since $df$ is nonsingular, $f$ is a 
local homeomorphism. Thus, since
$M_1\times M_2$ is compact,  it follows that $f$ is onto, and is 
therefore a covering map. But
$\S^{N-1}$ is simply connected (assuming
$n_1$, $n_2\neq 0$), so it admits no nontrivial coverings. This  means
that
$f$ has to be one-to-one, and hence a homeomorphism. Thus we obtain a 
contradiction, because, by basic
homology theory, a product of  manifolds cannot be homeomorphic to a
sphere.
\end{proof}

The next argument is more elementary:

\begin{proof}[Proof 2] 
Let $\conv M_1$ and $\conv M_2$ denote the 
convex hulls of $M_1$ and
$M_2$ respectively. First suppose that $\conv M_1=\conv M_2$.  Let 
$p_1$ be a point of $M_1$ in the
boundary of $\conv M_1$, and
$H\subset\R^N$ be a supporting hyperplane of $M_1$ through $p_1$. Then
$H$  also  supports $M_2$, and $H\cap
\conv(M_2)\neq\emptyset$; therefore,
$H$ must intersect
$M_2$ at some point $p_2$. Thus $T_{p_1} M_1$ and $T_{p_2} M_2$ both lie
in
$H\simeq\R^{N-1}$, and Lemma \ref{lem:A&G} yields that $N-1\geq 
n_1+n_2+1$, as desired.

So it remains to consider the case when $\conv M_1\neq\conv M_2$. 
Then one of our manifolds, say
$M_1$, must have a point, say $p_0$, which lies outside the 
convex hull of another.  In particular, $M_2$ must have a 
support hyperplane $H$ which separates
$M_2$ from $p_0$.  We may move $H$ parallel to itself and away from 
$\conv M_2$ until it becomes
tangent to $M_1$, say at a point
$p_1$. Since $\dim (T_{p_1} M_1)<\dim (H)$,  we may ``rotate"
$H$ around
$T_{p_1}M_1$ until it touches
$\conv M_2$ at some point $p_2$.  So once again we obtain a support 
hyperplane for
$M_1\cup M_2$ which intersects each submanifold, and we may  invoke 
Lemma \ref{lem:A&G} to complete the
proof.
\end{proof}

\section*{Acknowledgments}   
It is a  pleasure to acknowledge useful 
discussions with M. Farber, D. Fuchs,
M. Kossowski, K. Y. Lam, S. Lvovsky, B. Solomon, and S. Yuzvinsky. 
Part of this work was completed while the
first author was  at the University of South Carolina.  The  second
author acknowledges the hospitality of ETH in Zurich where part of this
work was done.

\end{document}